# The Cauchy Mean Value Theorem for Intuitionistic Fuzzy Functions

Efendi[1,*], Admi Nazra[1], Haripamyu[1], Mahdhivan Syafwan[1]

[1] Department of Mathematics and Data Sciences, Andalas University, West Sumatera, Indonesia

___

**Abstract**

In this paper, we first study the arithmetic properties of intuitionistic fuzzy number, the monotonicity of intuitionistic fuzzy function and the derivative of intuitionistic fuzzy functions and then we study the fundamental properties on monotone function. After that, we proposed a mean value theorem for differentiable function in intuitionistic fuzzy calculus, from which we derived a general version using Cauchy mean value theorem for intuitionistic fuzzy calculus. At last, we give some examples and we also show as a corollary that the mean value theorem is a special case of the Cauchy version of mean value theorem.

## 1. Introduction

Fuzzy sets, a groundbreaking concept introduced by Zadeh [1] in 1965, have left an indelible mark on various domains within modern society [2,3]. Rooted in their membership functions, these sets have proven indispensable for addressing the nuances of uncertainty and vagueness. However, their inherent limitation lies in their inability to fully encapsulate the multifaceted dimensions of uncertainty, which encompass not only membership but also non-membership, support, opposition, and hesitancy. To surmount this limitation, Atanassov [4] introduced intuitionistic fuzzy sets (IFS) in 1983, a profound extension of fuzzy sets that introduces both membership and non-membership functions. This expansion of the theoretical framework has found applications across domains like decision-making [5-9], pattern recognition [10-12], and medical diagnosis [13-14], underscoring its efficacy in modeling complex real-world uncertainties.

[*] Corresponding author. E-mail addresses: efendi@sci.unand.ac.id

In the past decades, the theory on IFSs has received great attention and been applied to decision making [2,18,19,20], pattern recognition [10,14,17], medical diagnosis [4,13] and so on. Inspired by the IFS, Xu and Yager [16] introduced the intuitionistic fuzzy number (IFN), which is an ordered non-negative number pair $(u, v)$, with $u + v \leq 1$, and gave some basic operations. Since then, Xu et.al published their works on intuitionistic fuzzy calculus [7-11]. Based on work of Xu et.al, then Yavuz, E. in [22], extend the concept of convergence of sequences in the domain of multiplicative calculus and its application to intuitionistic fuzzy numbers. In early 2023, Yavuz, E. then published A calculus for intuitionistic fuzzy values, that extended Xu et.al work and showed the relationship between intuitionistic fuzzy calculus and multiplicative that developed by Grossman M and Katz R.

The trajectory of intuitionistic fuzzy mathematics progresses with the introduction of intuitionistic fuzzy numbers (IFN) by Xu and Yager [16]. These IFNs, represented as ordered non-negative number pairs (u,v) with u ≤ v, provide a more refined means of capturing uncertainty and hesitancy. The crucial underpinnings of these numbers, encapsulating both membership and non-membership aspects, gain significance from the comprehensive operations associated with them [16]. This foundational work on IFNs broadens the mathematical toolkit available for navigating complex uncertainties, thereby serving as a stepping stone for further advancements.

Parallel to these developments, research has delved into the intricate realm of intuitionistic fuzzy calculus, which involves the exploration of derivatives, integrals, and their properties [17-21]. This exploration aligns with the overarching goal of refining the understanding and application of intuitionistic fuzzy constructs. Specifically, the study of derivative and differential operations of intuitionistic fuzzy numbers [17], along with fundamental properties of intuitionistic fuzzy calculus [18], refines the mathematical machinery that underpins this domain. These investigations offer a robust foundation for more complex analyses involving intuitionistic fuzzy functions and their derivatives.

Furthermore, relationships between various types of intuitionistic fuzzy integrals have been explored, shedding light on the nuances of aggregation and summation within uncertain environments [21,24]. The utilization of integrals based on Archimedean t-conorms and t-norms [22] showcases the adaptability and versatility of intuitionistic fuzzy mathematics in handling uncertainty. These investigations enhance the practical applicability of intuitionistic fuzzy mathematics, fostering more informed and comprehensive decision-making processes in real-world scenarios.

The journey through intuitionistic fuzzy mathematics encompasses a trajectory of innovation, extending from the foundational concepts of fuzzy sets to the sophisticated terrain of intuitionistic fuzzy numbers, calculus, and integrals. These advancements, rooted in the seminal works of pioneers like Zadeh [1] and Atanassov [4], continue to shape a mathematical landscape that resonates across diverse applications [2-14]. This ongoing exploration underscores the profound impact of intuitionistic fuzzy mathematics in addressing complex uncertainties and enriching the spectrum of mathematical tools available to researchers and practitioners alike.

This article stands at the forefront of extending the realm of intuitionistic fuzzy calculus, notably in the context of the mean value theorem. The generalization of the mean value theorem within the framework of intuitionistic fuzzy calculus constitutes a noteworthy contribution, bridging theoretical constructs with real-world applications. Building upon the foundation laid by previous works on intuitionistic fuzzy calculus [17-21], this article takes a significant stride by adapting and extending the well-established mean value theorem to the realm of intuitionistic fuzzy calculus. This extension opens new avenues for understanding the behavior of intuitionistic fuzzy functions, derivatives, and the intricate relationships between them.

By unraveling the generalized mean value theorem in intuitionistic fuzzy calculus, the article not only enhances our theoretical understanding of this mathematical domain but also equips us with a powerful tool for modeling complex uncertainties in practical scenarios. This extension is poised to have a profound impact on fields like decision-making, pattern recognition, and more, where the representation of multi-faceted uncertainties is pivotal. The article's contribution lies in its capacity to bridge the theoretical underpinnings of intuitionistic fuzzy mathematics with their application in real-world problem-solving, thereby strengthening the relevance and applicability of this evolving mathematical paradigm.

## 2. Preliminaries

In this section, we first introduce the concept of Zadeh's fuzzy set [1] and subsequently introduce intuitionistic fuzzy number, intuitionistic fuzzy function, intuitionistic fuzzy derivative [15,17], the mean value theorem and the Cauchy version of mean value theorem [22].

**Definition 2.1 [1].** Let $X$ be a fixed set, then $F = \{\langle x, u_F(x)\rangle \mid x \in X\}$ is called a fuzzy set, where $u_F$ is the membership function of $X$, $u_F : X \to [0,1]$, and $u_F(x)$ indicates the membership degree of the element $x$ to $F$, which is a single value belonging to the unit closed interval $[0,1]$.

It is noteworthy that the curly brace represents a set and the angular brace means that $x$ and $u_F(x)$ are taken as a whole. Besides, the nearer the value of the membership is to unity, the higher "grade of membership" of the element $x$ to $F$ is.

In 1983, Atanassov extended the fuzzy set to intuitionistic fuzzy set (IFS) [4] as defined below:

**Definition 2.2 [4].** Let $X$ be a fixed set, then an IFS is expressed as:
$A = \{\langle x, u_A(x), v_A(x)\rangle | x \in X\}$, each element of which is respectively depicted by the membership degree function: $u_A : X \to [0,1]$ and the non-membership degree function: $v_A : X \to [0,1]$ meeting the condition: $0 \leq u_A(x) + v_A(x) \leq 1, \forall x \in X$.

Besides, $u_A(x)$ and $v_A(x)$ represent, respectively, the membership degree and the non-membership degree of $x$ to $A$. When $1 - u_A(x) - v_A(x) = 0$, for all $x \in X$, then the IFS $A$ reduces to the fuzzy set. Similar to the membership, the nearer the value of the non-membership is to unity, the higher "grade of non-membership" of the element $x$ to $A$ is.

Xu and Yager [15] further defined the basic elements of an IFS as intuitionistic fuzzy number (IFN), which can be expressed by an ordered non-negative number pair $(u, v)$ meeting $u + v \leq 1$. Moreover, some basic operations of IFNs were also given as follows:

**Definition 2.3 [15, 17, 23].** Let $\alpha = (u_\alpha, v_\alpha), \beta = (u_\beta, v_\beta)$, be two IFNs, then

**(Addition)** $\alpha \oplus \beta = (u_\alpha + u_\beta - u_\alpha u_\beta, v_\alpha v_\beta)$

**(Subtraction)** $\alpha \ominus \beta = \begin{cases} \left(\frac{u_\alpha - u_\beta}{1 - u_\beta}, \frac{v_\alpha}{v_\beta}\right), & if\ 0 \leq \frac{v_\alpha}{v_\beta} \leq \frac{1 - u_\alpha}{1 - u_\beta} \leq 1 \\ (0, 1) & otherwise \end{cases}$

Xu and Yager [15] pointed out that the addition of IFNs satisfies the commutative law. In addition, Lei and Xu [20] further studied the changed region of the basic operations of IFNs, and denoted them in the $(u, v)$ plane. Firstly, we briefly introduce the changed region of the operations "$\oplus$" and "$\ominus$":

(1) Let $\alpha$ be a given IFN, then the region $S_\oplus(\alpha)$ represents $\alpha \oplus \varepsilon$ for any IFN $\varepsilon$, which is shown in Fig. 2.1;

(2) The region $S_\ominus(\alpha)$ represents $\alpha \ominus \varepsilon$ for any IFN $\varepsilon$, which is shown in Fig. 2.1.

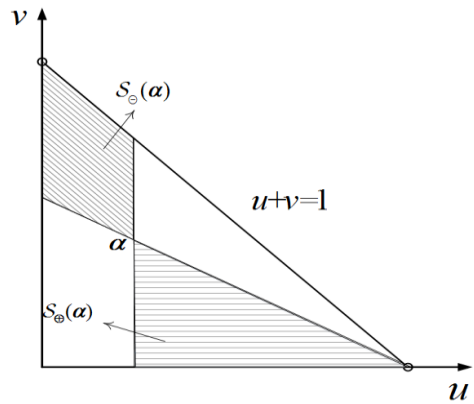

**Fig. 2.1 [17,20]. The regions** $S_\oplus(\alpha)$ **and** $S_\ominus(\alpha)$

In what follows, Xu et.al introduce the scalar-multiplication operation

**Definition 2.4 [15].** Let $\alpha = (u_\alpha, v_\alpha)$ be an IFN, and the parameter $\lambda$ be a real number meeting $\lambda > 0$. Then **(Scalar-multiplication)** $\lambda \alpha = (1 - (1 - u_\alpha)^\lambda, v_\alpha^\lambda)$

In order to understand the operation better, Xu et.al investigate some properties of the scalar-multiplication. For any given IFN $\alpha_0 = (u_0, v_0)$, and based on the mathematical expression of $\lambda \alpha_0$, Xu et.al have the following results [20]:

(1) $\lambda \alpha_0$ represents a concave curve in $(u, v)$ plane, and the concave curve is only decided by the function of $v(u) = v_0^{\frac{\ln(1-u)}{\ln(1-u_0)}}$ or $u(v) = 1 - (1 - u_0)^{\frac{\ln v}{\ln v_0}}$;

(2) When $\lambda\alpha_0 = (u,v)$, Xu et.al had shown that the parameter $\lambda$ if only $u_0 \neq 0$, $u_0 \neq 1$, $v_0 \neq 0$ and $v_0 \neq 1$;

(3) The value of $\lambda\alpha_0$ will change when $\lambda$ changes from zero to infinity.

Next, Xu et.al also introduced the properties of the function $v(u)$ in detail, and the properties about $u(v)$ is similar [20]:

(1) There is $v(u_0) = v_0$, which shows the situation when $\lambda = 1$;

(2) $v(1) = 0$ reveals that $\lambda\alpha_0 \to (1,0)$ when $\lambda \to +\infty$;

(3) $v(0) = 1$ corresponds to the situation $\lambda\alpha_0 \to (0,1)$ when $\lambda \to 0$;

(4) When $\lambda > 1$, there will be $\lambda\alpha_0 \in S_\oplus(\alpha_0)$ because $\lambda\alpha_0 = \alpha_0 \oplus (\lambda-1)\alpha_0$;

(5) When $0 < \lambda < 1$, there is $\lambda\alpha_0 \in S_\ominus(\alpha_0)$, due to that $\lambda\alpha_0 = \alpha_0 \ominus (1-\lambda)\alpha_0$.

Based on the above analyses, Xu et.al had shown that the graph $l_\alpha$ of $\lambda\alpha$, which is shown in Fig. 2.2 [20].

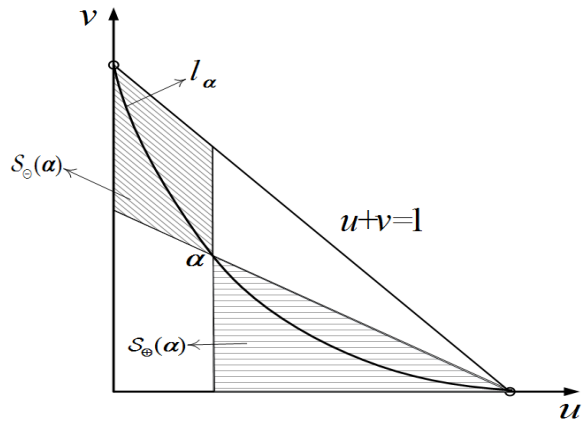

**Fig. 2.2 [20]. The graph of $\lambda\alpha$**

Especially, Xu et.al had shown that that addition, subtraction and scalar-multiplication are all closure operations from the changed region of the operations shown in Fig. 2.1 and Fig. 2.2.

In intuitionistic fuzzy calculus, an intuitionistic fuzzy function (IFFs) [19] has the form $\varphi(\alpha) = (f(\mu), g(v))$, which consists of two real number functions (i.e. $f(\mu)$ and $g(v)$) with the conditions: $0 \leq f(\mu), g(v) \leq 1$ and $0 \leq f(\mu) + g(v) \leq 1$.

Note: If there are no special instructions, the two functions are generally assumed to be continuous and derivable.

In the following, Xu et.al have introduced intuitionistic fuzzy calculus based on the subtraction derivative of intuitionistic fuzzy calculus,

**Definition 2.4** [19]. Let $X = (\mu, v)$, $\varphi(X) = (f(\mu), g(v))$ be an IFF, and $Y = X \oplus \Delta X \in S_\oplus(X)$, then

$$\frac{d^\oplus \varphi(X)}{dX} = \lim_{Y \to X} \left( \frac{\varphi(Y) \ominus \varphi(X)}{Y \ominus X} \right) = \lim_{\Delta X \to 0} \left( \frac{\varphi(X \oplus \Delta X) - \varphi(X)}{\Delta X} \right)$$

$$= \left( \frac{1-\mu}{1-f(\mu)} \frac{df(\mu)}{d\mu}, 1 - \frac{v}{g(v)} \frac{dg(v)}{dv} \right)$$

is defined as the addition derivative function of $\varphi(X)$ if the limit exist.

**Definition 2.5** [19]. Let $X = (\mu, v)$, $\varphi(X) = (f(\mu), g(v))$ be an IFF, and $Y = X \oplus \Delta X \in S_\oplus(X)$, then

$$\frac{d^\otimes \varphi(X)}{dX} = \lim_{Y \to X} \left( \frac{\varphi(Y)}{\varphi(X)} \ominus \frac{Y}{X} \right) = \lim_{\Delta X \to 0} \left( \frac{\varphi(X \otimes \Delta X)}{\varphi(X)} \ominus \Delta X \right)$$

$$= \left( 1 - \frac{\mu}{f(\mu)} \frac{df(\mu)}{d\mu}, \frac{1-v}{1-g(v)} \frac{dg(v)}{dv} \right)$$

is defined as the multiplication derivative function of $\varphi(X)$ if the limit exist.

If there is an IFN $\gamma$, which satisfies $\alpha \oplus \gamma = \beta$, then $\alpha$ is less than or equal to $\beta$, denoted by $\alpha \trianglelefteq \beta$. If there is an IFN $\gamma$, which satisfies $\alpha \oplus \gamma = \beta$ and $\gamma \neq (0,1) = E$, then $\alpha$ is less than $\beta$, denoted by $\alpha \triangleleft \beta$.

**Definition 2.6**. If $\alpha \trianglelefteq \beta$, then $\varphi(\alpha) \trianglelefteq \varphi(\beta)$, in this case $\varphi(\alpha)$ a monotone increasing IFF.

In fact, the intuitionistic fuzzy calculus is mainly on the basis of the monotone increasing IFFs.

**Definition 2.7.** Let $\varphi(\alpha)$ be an IFF that exists the subtraction derivative function $\frac{d\varphi(\alpha)}{d\alpha}$. If noting $\Delta\alpha = \alpha' \ominus \alpha$, then $d\varphi(\alpha) = \frac{d\varphi(\alpha)}{d\alpha} \otimes \Delta\alpha$ is called the subtraction differential of $\varphi(\alpha)$. Since $d\alpha = \Delta\alpha$, then the subtraction differential can be represented as $d\varphi(\alpha) = \frac{d\varphi(\alpha)}{d\alpha} \otimes d\alpha$.

Based on the differential of IFFs above, the following theorem will establish the relationship between the function change value ($\Delta\varphi(\alpha)$) and the differential ($d\varphi(\alpha)$).

### 2.1 Mean Value Theorem

In Mathematics, the mean value theorem (or Lagrange theorem) states, roughly, that for a given planar arc between two endpoints, there is at least one point at which the tangent to the arc is parallel to the secant through its endpoints. It is one of the most important results in real analysis. This theorem is used to prove statements about a function on an interval starting from local hypotheses about derivatives at points of the interval.

More precisely, the theorem states that if the functions $f$ and $g$ are both continuous on the closed interval $[a, b]$ and differentiable on the open interval $(a, b)$, then there exists some $c \in (a, b)$ such that:

$$f'(c) = \frac{f(b) - f(a)}{b - a}$$

Theorem 2.1 below stated an Intuitionistic fuzzy properties on monotone function proposed by Xu et.al, by which this paper is mainly inspired. In the next section, we will slightly modify Theorem 2.1 into mean value theorem and then generalize it into cauchy mean value theorem.

**Theorem 2.1** [19]. [Zeshui Xu Monotone Theorem] Let $\varphi(X) = (f(\mu), g(v))$ be a monotone increasing IFF of $X \trianglelefteq Y$. If $\varphi(X)$ is differentiable then:

$$\varphi(Y) \ominus \varphi(X) \approx \frac{d\varphi(X)}{dX} (Y \ominus X)$$

**Proof:**

If $X \trianglelefteq Y$, $X$ and $Y$ are expressed by $(\mu, v)$ and $(\mu', v')$ respectively, then

$$\Delta X = Y \ominus X = (\mu', v') \ominus (\mu, v) = \left(\frac{\mu' - \mu}{1 - \mu}, \frac{v'}{v}\right)$$

and also,

$$\frac{d\varphi(X)}{dX} = \left(\frac{1 - \mu}{1 - f(\mu)} \frac{df(\mu)}{d\mu}, 1 - \frac{v}{g(v)} \frac{dg(v)}{dv}\right)$$

Hence, we have:

$$\frac{d\varphi(X)}{dX} \otimes \Delta X$$

$$= \left(\frac{1 - \mu}{1 - f(\mu)} \frac{df(\mu)}{d\mu}, 1 - \frac{v}{g(v)} \frac{dg(v)}{dv}\right) \otimes \left(\frac{\mu' - \mu}{1 - \mu}, \frac{v'}{v}\right)$$

$$= \left(\frac{1 - \mu}{1 - f(\mu)} \frac{df(\mu)}{d\mu} \frac{\mu' - \mu}{1 - \mu}, 1 - \frac{v}{g(v)} \frac{dg(v)}{dv} + \frac{v'}{v} - \left(1 - \frac{v}{g(v)} \frac{dg(v)}{dv}\right) \frac{v'}{v}\right)$$

$$= \left(\frac{\mu' - \mu}{1 - f(\mu)} \frac{df(\mu)}{d\mu}, 1 - \frac{v}{g(v)} \frac{dg(v)}{dv} + \frac{v'}{v} - \left(\frac{v'}{v} - \frac{v}{g(v)} \frac{dg(v)}{dv} \frac{v'}{v}\right)\right)$$

$$= \left(\frac{\mu' - \mu}{1 - f(\mu)} \frac{df(\mu)}{d\mu}, 1 - \frac{v}{g(v)} \frac{dg(v)}{dv} + \left(\frac{v'}{g(v)} \frac{dg(v)}{dv}\right)\right)$$

$$= \left(\frac{\mu' - \mu}{1 - f(\mu)} \frac{df(\mu)}{d\mu}, 1 - \frac{v - v'}{g(v)} \frac{dg(v)}{dv}\right)$$

Therefore:

$$\varphi(X) \oplus \frac{d\varphi(X)}{dX} \otimes \Delta X$$

$$= (f(\mu), g(v)) \oplus \left(\frac{\mu' - \mu}{1 - f(\mu)} \frac{df(\mu)}{d\mu}, 1 - \frac{v - v'}{g(v)} \frac{dg(v)}{dv}\right)$$

$$= \left(f(\mu) + \frac{\mu' - \mu}{1 - f(\mu)} \frac{df(\mu)}{d\mu} - f(\mu) \left(\frac{\mu' - \mu}{1 - f(\mu)} \frac{df(\mu)}{d\mu}\right), g(v) \left(1 - \frac{v - v'}{g(v)} \frac{dg(v)}{dv}\right)\right)$$

$$= \left(f(\mu) + (1 - f(\mu))\left(\frac{\mu' - \mu}{1 - f(\mu)}\frac{df(\mu)}{d\mu}\right), g(v) + (v' - v)\frac{dg(v)}{dv}\right)$$

$$= \left(f(\mu) + (\mu' - \mu)\frac{df(\mu)}{d\mu}, g(v) + (v' - v)\frac{dg(v)}{dv}\right)$$

$$= \left(f(\mu) + (\mu' - \mu)\left(\frac{f(\mu') - f(\mu)}{\mu' - \mu} + o\right), g(v) + (v' - v)\left(\frac{g(v') - g(v)}{v' - v} + o\right)\right)$$

$$= \left(f(\mu') + o(\mu' - \mu), g(v') + o(v' - v)\right)$$

$$\approx \varphi(Y)$$

Hence,

$$\varphi(X) \oplus \frac{d\varphi(X)}{dX} \otimes \Delta X \approx \varphi(Y)$$

Or

$$\varphi(Y) \ominus \varphi(X) \approx \frac{d\varphi(X)}{dX}(Y \ominus X)$$

QED.

## 2.2 Cauchy Mean value theorem

The mean value theorem is the special case of Cauchy's mean value theorem when $g(t) = t$. Cauchy's mean value theorem, also known as the extended mean value theorem, is a generalization of the mean value theorem. It states: *if the functions $f$ and $g$ are both continuous on the closed interval $[a, b]$ and differentiable on the open interval $(a, b)$, then there exists some $c \in (a, b)$ such that $(f(b) - f(a))g'(c) = (g(b) - g(a))f'(c)$.*

Of course, if $g(b) \neq g(a)$ and $g'(c) \neq 0$ this is equivalent to: $\frac{f'(c)}{g'(c)} = \frac{f(b)-f(a)}{g(b)-g(a)}$. Geometrically, this means that there is some tangent to the graph of the curve $t \to (f(t), g(t))$ which is parallel to the line defined by the points $(f(a), g(a))$ and $(f(b), g(b))$. However, Cauchy's theorem does not claim the existence of such a tangent in all cases where $(f(a), g(a))$ and $(f(b), g(b))$ are distinct points, since it might be satisfied only for some value c with $f'(c) = g'(c) = 0$, in other words a value for which the mentioned curve is stationary; in such points no tangent to the curve is likely to be defined at all. An example of this situation is the curve given by $t \to (t^3, 1 - t^2)$. which on the interval $[-1,1]$ goes from the point $(-1,0)$ to $(1,0)$, yet never has a horizontal tangent; however it has a stationary point (in fact a cusp) at $t = 0$. Suppose $g(a) \neq g(b)$. Define $h(x) = f(x) - rg(x)$, where $r = \frac{f(b)-f(a)}{g(b)-g(a)}$ is fixed. Since $f$ and $g$ are continuous on $[a, b]$ and differentiable on $(a, b)$, the same is true for $h$. All in all, $h$ satisfies the conditions of Rolle's theorem: consequently, there is some c in $(a, b)$ for which $h'(c) = 0$. Now using the definition of h we have:

$$0 = h'(c) = f'(c) - rg'(c) = f'(c) - \frac{f(b) - f(a)}{g(b) - g(a)} g'(c)$$

Therefore:

$$f'(c) = \frac{f(b) - f(a)}{g(b) - g(a)} g'(c)$$

which implies the result.

## 3. Main Results

In this section, we will modify theorem 2.1 into MVT and then we will generalize it into Cauchy version. We will show a proof of Cauchy's mean value theorem. The proof of mean value theorem (MVT) and Cauchy's mean value theorem (CMVT) is based on the same idea as the proof of the Xu theorem for monotone function (XMT). And the last we derive some corollaries from CMVT. The most important corollary is that MVT is a special case of CMVT. We will divide this section into two subsections, Addition CMVT and Multiplication CMVT.

### 3.1 Addition CMVT

**Theorem 3.1** [MVT} Let $\varphi(X) = (f(\mu), g(v))$ be a monotone increasing IFF of $X \trianglelefteq Y$. If $\varphi(X)$ is differentiable, then:

$$\varphi(Y) \ominus \varphi(X) = \frac{d^{\oplus}\varphi(X)}{dX}(Y \ominus X)$$

**Proof:**

If $X \trianglelefteq Y$, $X$ and $Y$ are expressed by $(\mu, v)$ and $(\mu', v')$ respectively, then

$$\Delta X = Y \ominus X = (\mu', v') \ominus (\mu, v) = \left(\frac{\mu' - \mu}{1 - \mu}, \frac{v'}{v}\right)$$

and also,

$$\frac{d^{\oplus}\varphi(X)}{dX} = \left(\frac{1 - \mu}{1 - f(\mu)} \frac{df(\mu)}{d\mu}, 1 - \frac{v}{g(v)} \frac{dg(v)}{dv}\right)$$

Hence, we have:

$$\varphi(Y) \ominus \varphi(X) = \big(f(\mu'), g(v')\big) \ominus \big(f(\mu), g(v)\big)$$
$$= \left(\frac{f(\mu') - f(\mu)}{1 - f(\mu)}, \frac{g(v')}{g(v)}\right)$$

And,

$$\frac{d^{\oplus}\varphi(X)}{dX} \otimes \Delta X = \left(\frac{1-\mu}{1-f(\mu)}\frac{df(\mu)}{d\mu}, 1 - \frac{v}{g(v)}\frac{dg(v)}{dv}\right) \otimes \left(\frac{\mu'-\mu}{1-\mu}, \frac{v'}{v}\right)$$

$$= \left(\frac{1-\mu}{1-f(\mu)}\frac{df(\mu)}{d\mu}\frac{\mu'-\mu}{1-\mu}, 1 - \frac{v}{g(v)}\frac{dg(v)}{dv} + \frac{v'}{v} - \left(1 - \frac{v}{g(v)}\frac{dg(v)}{dv}\right)\frac{v'}{v}\right)$$

$$= \left(\frac{\mu'-\mu}{1-f(\mu)}\frac{df(\mu)}{d\mu}, 1 - \frac{v}{g(v)}\frac{dg(v)}{dv} + \frac{v'}{v} - \left(\frac{v'}{v} - \frac{v}{g(v)}\frac{dg(v)}{dv}\frac{v'}{v}\right)\right)$$

$$= \left(\frac{\mu'-\mu}{1-f(\mu)}\frac{df(\mu)}{d\mu}, 1 - \frac{v}{g(v)}\frac{dg(v)}{dv} + \left(\frac{v'}{g(v)}\frac{dg(v)}{dv}\right)\right)$$

$$= \left(\frac{\mu'-\mu}{1-f(\mu)}\frac{df(\mu)}{d\mu}, 1 - \frac{v-v'}{g(v)}\frac{dg(v)}{dv}\right)$$

$$= \left(\frac{\mu' - \mu}{1 - f(\mu)} \frac{f(\mu') - f(\mu)}{\mu' - \mu}, 1 - \frac{v - v'}{g(v)} \frac{g(v') - g(v)}{v' - v}\right)$$

$$= \left(\frac{\mu' - \mu}{1 - f(\mu)} \frac{f(\mu') - f(\mu)}{\mu' - \mu}, 1 + \frac{g(v') - g(v)}{g(v)}\right)$$

$$= \left(\frac{\mu' - \mu}{1 - f(\mu)} \frac{f(\mu') - f(\mu)}{\mu' - \mu}, 1 + \frac{g(v')}{g(v)} - 1\right) = \left(\frac{f(\mu') - f(\mu)}{1 - f(\mu)}, \frac{g(v')}{g(v)}\right)$$

$$= (f(\mu'), g(v')) \ominus (f(\mu), g(v)) = \varphi(Y) \ominus \varphi(X)$$

Therefore:

$$\varphi(Y) \ominus \varphi(X) = \frac{d^\oplus \varphi(X)}{dX}(Y \ominus X)$$

<div style="text-align: right;">QED.</div>

**Example 1:**

Prove that MVT can be written as: $\frac{d^\oplus \varphi(X_0)}{dX} = \frac{\varphi(Y) \ominus \varphi(X)}{(Y \ominus X)}$

**Proof**:

$$\frac{\varphi(Y) \ominus \varphi(X)}{(Y \ominus X)} = \frac{(f(\mu'), g(v')) \ominus (f(\mu), g(v))}{(Y \ominus X)} = \frac{\left(\frac{f(\mu') - f(\mu)}{1 - f(\mu)}, \frac{g(v')}{g(v)}\right)}{\left(\frac{\mu' - \mu}{1 - \mu}, \frac{v'}{v}\right)}$$

$$= \frac{\left(\frac{f(\mu') - f(\mu)}{1 - f(\mu)}, \frac{g(v')}{g(v)}\right)}{\left(\frac{\mu' - \mu}{1 - \mu}, \frac{v'}{v}\right)} = \left(\frac{\frac{f(\mu') - f(\mu)}{1 - f(\mu)}}{\frac{\mu' - \mu}{1 - \mu}}, \frac{\frac{g(v')}{g(v)} - \frac{v'}{v}}{1 - \frac{v'}{v}}\right)$$

$$= \left(\frac{f(\mu') - f(\mu)}{1 - f(\mu)} \cdot \frac{1 - \mu}{\mu' - \mu}, \frac{\frac{vg(v')}{vg(v)} - \frac{v'g(v)}{vg(v)}}{\frac{v - v'}{v}}\right)$$

$$= \left(\frac{f(\mu') - f(\mu)}{1 - f(\mu)} \cdot \frac{1 - \mu}{\mu' - \mu}, \frac{\frac{1}{g(v)}(vg(v') - v'g(v))}{v - v'}\right)$$

$$= \left(\frac{f(\mu') - f(\mu)}{1 - f(\mu)} \cdot \frac{1 - \mu}{\mu' - \mu}, \frac{1}{g(v)} \cdot \frac{(vg(v') - vg(v) + vg(v) - v'g(v))}{v - v'}\right)$$

$$= \left(\frac{f(\mu') - f(\mu)}{1 - f(\mu)} \cdot \frac{1 - \mu}{\mu' - \mu}, \frac{1}{g(v)}\left(-v\frac{(g(v') - g(v))}{v' - v} + \frac{(v - v')}{v - v'} g(v)\right)\right)$$

$$= \left(\frac{1 - \mu}{1 - f(\mu)} \cdot \frac{f(\mu') - f(\mu)}{\mu' - \mu}, \frac{1}{g(v)}\left(-v\frac{(g(v') - g(v))}{v' - v} + g(v)\right)\right)$$

$$= \left(\frac{1 - \mu}{1 - f(\mu)} \cdot \frac{f(\mu') - f(\mu)}{\mu' - \mu}, 1 - \frac{v}{g(v)} \cdot \frac{g(v') - g(v)}{v' - v}\right)$$

$$= \left(\frac{1 - \mu}{1 - f(\mu)} \cdot \frac{df(\mu_0)}{d\mu}, 1 - \frac{v}{g(v)} \cdot \frac{dg(v_0)}{dv}\right) = \frac{d^{\oplus}\varphi(X_0)}{dX}$$

**Example 2:**

Find $X_0$ where $X \triangleleft X_0 \triangleleft Y$ and $\dfrac{d^{\oplus}\varphi(X_0)}{dX} = \dfrac{\varphi(Y) \ominus \varphi(X)}{(Y \ominus X)}$

**Solution**

Let:

$$\varphi(X) = X^2 = (\mu, v)^2 = (\mu^2, 1 - (1 - v)^2) = (f(\mu), g(v))$$

$$X_1 = (\mu_1, v_1) = (0.1, 0.7)$$
$$Y_1 = (\mu_1', v_1') = (0.6, 0.3)$$
$$\varphi(X_1) = X_1^2 = (0.1, 0.7)^2 = (0.01, 0.91) = (f(\mu), g(v))$$
$$\varphi(Y_1) = Y_1^2 = (0.6, 0.3)^2 = (0.36, 0.51) = (f(\mu'), g(v'))$$
$$\varphi(Y) \ominus \varphi(X) = \left( \frac{1-\mu}{1-f(\mu)} \frac{f(\mu') - f(\mu)}{\mu' - \mu}, 1 - \frac{v}{g(v)} \frac{g(v') - g(v)}{v' - v} \right)$$
$$\varphi(Y) \ominus \varphi(X) = \left( \frac{1-0.1}{1-0.01} \frac{0.36 - 0.01}{0.6 - 0.1}, 1 - \frac{0.7}{0.91} \frac{0.51 - 0.91}{0.3 - 0.7} \right)$$
$$\frac{d^\oplus \varphi(X_0)}{dX} = \left( \frac{1-\mu}{1-f(\mu)} \frac{df(\mu_0)}{d\mu}, 1 - \frac{v}{g(v)} \frac{dg(v_0)}{dv} \right)$$
$$\frac{df(\mu_0)}{d\mu} = \frac{0.36 - 0.01}{0.6 - 0.1}$$
$$2\mu_0 = 0.7$$
$$\mu_0 = 0.35$$

$$\frac{dg(v_0)}{dv} = \frac{0.51 - 0.91}{0.3 - 0.7}$$
$$2(1 - v_0) = 1$$
$$v_0 = 0.5$$

So that, $X_0 = (\mu_0, v_0) = (0.35, 0.5)$

**Example 3:**

Find $X_0$ where $X \triangleleft X_0 \triangleleft Y$ and $\dfrac{d^\oplus \varphi(X_0)}{dX} = \dfrac{\varphi(Y) \ominus \varphi(X)}{(Y \ominus X)}$

**Solution**

Let:
$$\varphi(X) = X^3 = (\mu, v)^3 = (\mu^3, 1 - (1-v)^3) = (f(\mu), g(v))$$

$$X_1 = (\mu_1, v_1) = (0.1, 0.7)$$

$$Y_1 = (\mu'_1, v'_1) = (0.6, 0.3)$$

$$\varphi(X_1) = X_1^{\ 3} = (0.1, 0.7)^3 = (0.001, 0.973) = (f(\mu), g(v))$$

$$\varphi(Y_1) = Y_1^{\ 3} = (0.6, 0.3)^3 = (0.216, 0.657) = (f(\mu'), g(v'))$$

$$\varphi(Y) \ominus \varphi(X) = \left( \frac{1-\mu}{1-f(\mu)} \frac{f(\mu') - f(\mu)}{\mu' - \mu}, 1 - \frac{v}{g(v)} \frac{g(v') - g(v)}{v' - v} \right)$$

$$\varphi(Y_1) \ominus \varphi(X_1) = \left( \frac{1-0.1}{1-0.001} \frac{0.216 - 0.001}{0.6 - 0.1}, 1 - \frac{0.7}{0.973} \frac{0.657 - 0.973}{0.3 - 0.7} \right)$$

$$\frac{d^{\oplus}\varphi(X_0)}{dX} = \left( \frac{1-\mu}{1-f(\mu)} \frac{df(\mu_0)}{d\mu}, 1 - \frac{v}{g(v)} \frac{dg(v_0)}{dv} \right)$$

$$\frac{df(\mu_0)}{d\mu} = \frac{0.216 - 0.001}{0.6 - 0.1}$$

$$3\mu_0^{\ 2} = 0.43$$

$$\mu_0 = 0.378593889$$

$$\frac{dg(v_0)}{dv} = \frac{0.657 - 0.973}{0.3 - 0.7}$$

$$3(1 - v_0)^2 = 0{,}79$$

$$v_0 = 0.486839889$$

So that, $X_0 = (\mu_0, v_0) = (0.378593889, 0.486839889)$

Based on theorem 3.1 we generalize MVT into CMVT, as follows:

**Theorem 3.2**

Let $\varphi(X) = (f(\mu), g(v))$ and $\gamma(X) = (f_1(\mu), g_1(v))$ be a monotone increasing IFF of $X \unlhd Y$. If $\varphi(X)$ and $\gamma(X)$ are differentiable, then

$$[\varphi(Y) \ominus \varphi(X)] \otimes \frac{d^{\oplus}\gamma(X)}{dX} = \frac{d^{\oplus}\varphi(X)}{dX} \otimes (\gamma(Y) \ominus \gamma(X))$$

***Proof:***

Let:
$$X = (\mu, v)$$
$$Y = (\mu', v')$$
$$\varphi(X) = (f(\mu), g(v))$$
$$\varphi(Y) = (f(\mu'), g(v'))$$
$$\gamma(X) = (f_1(\mu), g_1(v))$$
$$\gamma(Y) = (f_1(\mu'), g_1(v'))$$
$$\frac{d^\oplus \varphi(X)}{dX} = \left(\frac{1-\mu}{1-f(\mu)}\frac{df(\mu)}{d\mu}, 1 - \frac{v}{g(v)}\frac{dg(v)}{dv}\right)$$
$$\frac{d^\oplus \gamma(X)}{dX} = \left(\frac{1-\mu}{1-f_1(\mu)}\frac{df_1(\mu)}{d\mu}, 1 - \frac{v}{g_1(v)}\frac{dg_1(v)}{dv}\right)$$

Then we have:
$$\varphi(Y) \ominus \varphi(X)$$
$$= (f(\mu'), g(v')) \ominus (f(\mu), g(v))$$
$$= \left(\frac{f(\mu') - f(\mu)}{1 - f(\mu)}, \frac{g(v')}{g(v)}\right)$$

And also:
$$\gamma(Y) \ominus \gamma(X)$$
$$= (f_1(\mu'), g_1(v')) \ominus (f_1(\mu), g_1(v))$$
$$= \left(\frac{f_1(\mu') - f_1(\mu)}{1 - f_1(\mu)}, \frac{g_1(v')}{g_1(v)}\right)$$

So that, the left hand side is:
$$[\varphi(Y) \ominus \varphi(X)] \otimes \frac{d^\oplus \gamma(X)}{dX}$$
$$= \left(\frac{f(\mu') - f(\mu)}{1 - f(\mu)}, \frac{g(v')}{g(v)}\right) \otimes \left(\frac{1-\mu}{1-f_1(\mu)}\frac{df_1(\mu)}{d\mu}, 1 - \frac{v}{g_1(v)}\frac{dg_1(v)}{dv}\right)$$

$$= \left( \frac{f(\mu') - f(\mu)}{1 - f(\mu)} \left( \frac{1 - \mu}{1 - f_1(\mu)} \frac{df_1(\mu)}{d\mu} \right) \right),$$

$$\frac{g(v')}{g(v)} + \left( 1 - \frac{v}{g_1(v)} \frac{dg_1(v)}{dv} \right) - \frac{g(v')}{g(v)} \left( 1 - \frac{v}{g_1(v)} \frac{dg_1(v)}{dv} \right) \right)$$

$$= \left( \frac{f(\mu') - f(\mu)}{1 - f(\mu)} \left( \frac{1 - \mu}{1 - f_1(\mu)} \frac{df_1(\mu)}{d\mu} \right) \right),$$

$$\frac{g(v')}{g(v)} + \left( 1 - \frac{v}{g_1(v)} \frac{dg_1(v)}{dv} \right) - \frac{g(v')}{g(v)} + \frac{g(v')}{g(v)} \left( \frac{v}{g_1(v)} \frac{dg_1(v)}{dv} \right) \right)$$

$$= \left( \frac{f(\mu') - f(\mu)}{1 - f(\mu)} \left( \frac{1 - \mu}{1 - f_1(\mu)} \frac{df_1(\mu)}{d\mu} \right), \left( 1 - \frac{v}{g_1(v)} \frac{dg_1(v)}{dv} \right) + \frac{g(v')}{g(v)} \left( \frac{v}{g_1(v)} \frac{dg_1(v)}{dv} \right) \right)$$

$$= \left( \frac{f(\mu') - f(\mu)}{1 - f(\mu)} \left( \frac{1 - \mu}{1 - f_1(\mu)} \frac{df_1(\mu)}{d\mu} \right), \left( 1 + \left( \frac{g(v')}{g(v)} - 1 \right) \frac{v}{g_1(v)} \frac{dg_1(v)}{dv} \right) \right)$$

$$= \left( \frac{f(\mu') - f(\mu)}{1 - f(\mu)} \left( \frac{1 - \mu}{1 - f_1(\mu)} \frac{df_1(\mu)}{d\mu} \right), \left( 1 + \left( \frac{g(v') - g(v)}{g(v)} \right) \frac{v}{g_1(v)} \frac{dg_1(v)}{dv} \right) \right)$$

$$= \left( \frac{f(\mu') - f(\mu)}{1 - f(\mu)} \left( \frac{1 - \mu}{1 - f_1(\mu)} \left( \frac{f_1(\mu') - f_1(\mu)}{\mu' - \mu} \right) \right), \left( 1 + \left( \frac{g(v') - g(v)}{g(v)} \right) \frac{v}{g_1(v)} \left( \frac{g_1(v') - g_1(v)}{v' - v} \right) \right) \right)$$

$$= \left( \frac{f(\mu') - f(\mu)}{1 - f(\mu)} \left( \frac{1 - \mu}{1 - f_1(\mu)} \left( \frac{f_1(\mu') - f_1(\mu)}{\mu' - \mu} \right) \right), \left( 1 + \left( \frac{g(v') - g(v)}{g(v)} \right) \frac{v}{g_1(v)} \left( \frac{g_1(v') - g_1(v)}{v' - v} \right) \right) \right)$$

$$= \left(\left(\left(\frac{f(\mu') - f(\mu)}{1 - f(\mu)}\right)\left(\frac{f_1(\mu') - f_1(\mu)}{1 - f_1(\mu)}\right)\left(\frac{1 - \mu}{\mu' - \mu}\right), 1\right.\right.$$

$$\left.\left.+ \left(\frac{g(v') - g(v)}{g(v)}\right)\left(\frac{g_1(v') - g_1(v)}{g_1(v)}\right)\left(\frac{v}{v' - v}\right)\right)\right)$$

The right hand side is:

$$[(\gamma(Y) \ominus \gamma(X))] \otimes \frac{d^\oplus \varphi(X)}{dX}$$

$$= \left(\frac{f_1(\mu') - f_1(\mu)}{1 - f_1(\mu)}, \frac{g_1(v')}{g_1(v)}\right) \otimes \left(\frac{1 - \mu}{1 - f(\mu)}\frac{df(\mu)}{d\mu}, 1 - \frac{v}{g(v)}\frac{dg(v)}{dv}\right)$$

$$= \left(\frac{f_1(\mu') - f_1(\mu)}{1 - f_1(\mu)}\left(\frac{1 - \mu}{1 - f(\mu)}\frac{df(\mu)}{d\mu}\right),\right.$$

$$\left.\frac{g_1(v')}{g_1(v)} + \left(1 - \frac{v}{g(v)}\frac{dg(v)}{dv}\right) - \frac{g_1(v')}{g_1(v)} + \left(1 - \frac{v}{g(v)}\frac{dg(v)}{dv}\right)\right)$$

$$= \left(\frac{f_1(\mu') - f_1(\mu)}{1 - f_1(\mu)}\left(\frac{1 - \mu}{1 - f(\mu)}\frac{df(\mu)}{d\mu}\right),\right.$$

$$\left.\frac{g_1(v')}{g_1(v)} + \left(1 - \frac{v}{g(v)}\frac{dg(v)}{dv}\right) - \frac{g_1(v')}{g_1(v)} + \frac{g_1(v')}{g_1(v)}\left(\frac{v}{g(v)}\frac{dg(v)}{dv}\right)\right)$$

$$= \left(\frac{f_1(\mu') - f_1(\mu)}{1 - f_1(\mu)}\left(\frac{1 - \mu}{1 - f(\mu)}\frac{df(\mu)}{d\mu}\right), \left(1 - \frac{v}{g(v)}\frac{dg(v)}{dv}\right) + \frac{g_1(v')}{g_1(v)}\left(\frac{v}{g(v)}\frac{dg(v)}{dv}\right)\right)$$

$$= \left(\frac{f_1(\mu') - f_1(\mu)}{1 - f_1(\mu)}\left(\frac{1 - \mu}{1 - f(\mu)}\frac{df(\mu)}{d\mu}\right), 1 + \left(\frac{g_1(v')}{g_1(v)} - 1\right)\left(\frac{v}{g(v)}\frac{dg(v)}{dv}\right)\right)$$

$$= \left(\frac{f_1(\mu') - f_1(\mu)}{1 - f_1(\mu)}\left(\frac{1 - \mu}{1 - f(\mu)}\frac{df(\mu)}{d\mu}\right), 1 + \left(\frac{g_1(v') - g_1(v)}{g_1(v)}\right)\left(\frac{v}{g(v)}\frac{dg(v)}{dv}\right)\right)$$

$$= \left(\frac{f_1(\mu') - f_1(\mu)}{1 - f_1(\mu)}\left(\frac{1 - \mu}{1 - f(\mu)}\left(\frac{f(\mu') - f(\mu)}{\mu' - \mu}\right)\right), 1\right.$$

$$\left.+ \left(\frac{g_1(v') - g_1(v)}{g_1(v)}\right)\left(\frac{v}{g(v)}\left(\frac{g(v') - g(v)}{v' - v}\right)\right)\right)$$

$$= \left( \frac{f_1(\mu') - f_1(\mu)}{1 - f_1(\mu)} \left( \frac{1 - \mu}{1 - f(\mu)} \left( \frac{f(\mu') - f(\mu)}{\mu' - \mu} \right) \right), 1 \right.$$

$$\left. + \left( \frac{g_1(v') - g_1(v)}{g_1(v)} \right) \left( \frac{v}{g(v)} \left( \frac{g(v') - g(v)}{v' - v} \right) \right) \right)$$

$$= \left( \left( \frac{f_1(\mu') - f_1(\mu)}{1 - f_1(\mu)} \right) \left( \frac{f(\mu') - f(\mu)}{1 - f(\mu)} \right) \left( \frac{1 - \mu}{\mu' - \mu} \right), 1 \right.$$

$$\left. + \left( \frac{g_1(v') - g_1(v)}{g_1(v)} \right) \left( \frac{g(v') - g(v)}{g(v)} \right) \left( \frac{v}{v' - v} \right) \right)$$

Now, the left hand side and the right hand side are equal.

**QED**

For an illustration, we give two examples:

**Example 1:**

Let

$$X = (\mu, v)$$
$$Y = (\mu', v')$$
$$\varphi(X) = X^2 = (\mu, v)^2 = (\mu^2, 1 - (1 - v)^2) = (f(\mu), g(v))$$
$$\varphi(Y) = Y^2 = (\mu', v')^2 = (\mu'^2, 1 - (1 - v')^2) = (f(\mu'), g(v'))$$
$$\gamma(X) = X^3 = (\mu, v)^3 = (\mu^3, 1 - (1 - v)^3) = (f_1(\mu), g_1(v))$$
$$\gamma(Y) = Y^3 = (\mu', v')^3 = (\mu'^3, 1 - (1 - v')^3) = (f_1(\mu'), g_1(v'))$$

So that, the left hand side:

$$[\varphi(Y) \ominus \varphi(X)] \otimes \frac{d^\oplus \gamma(X)}{dX}$$

$$= \left( \left( \frac{f(\mu') - f(\mu)}{1 - f(\mu)} \right) \left( \frac{f_1(\mu') - f_1(\mu)}{1 - f_1(\mu)} \right) \left( \frac{1 - \mu}{\mu' - \mu} \right), 1 \right.$$

$$\left. + \left( \frac{g(v') - g(v)}{g(v)} \right) \left( \frac{g_1(v') - g_1(v)}{g_1(v)} \right) \left( \frac{v}{v' - v} \right) \right)$$

$$= \left(\frac{\mu'^3 - \mu^3}{1 - \mu^3}\left(\frac{\mu'^2 - \mu^2}{1 - \mu^2}\right)\left(\frac{1-\mu}{\mu'-\mu}\right), 1\right.$$

$$\left. + \left(\frac{(1-(1-v')^3) - (1-(1-v)^3)}{(1-(1-v)^3)}\right)\left(\frac{(1-(1-v')^2) - ((1-(1-v)^2))}{(1-(1-v)^2)}\right)\left(\frac{v}{v'-v}\right)\right)$$

The right hand side:

$$(\gamma(Y) \ominus \gamma(X)) \otimes \frac{d^\oplus \varphi(X)}{dX}$$

$$= \left(\left(\frac{f_1(\mu') - f_1(\mu)}{1 - f_1(\mu)}\right)\left(\frac{f(\mu') - f(\mu)}{1 - f(\mu)}\right)\left(\frac{1-\mu}{\mu'-\mu}\right), 1\right.$$

$$\left. + \left(\frac{g_1(v') - g_1(v)}{g_1(v)}\right)\left(\frac{g(v') - g(v)}{g(v)}\right)\left(\frac{v}{v'-v}\right)\right)$$

$$= \left(\left(\frac{\mu'^2 - \mu^2}{1 - \mu^2}\right)\left(\frac{\mu'^3 - \mu^3}{1 - \mu^3}\right)\left(\frac{1-\mu}{\mu'-\mu}\right), 1\right.$$

$$\left. + \left(\frac{(1-(1-v')^2) - ((1-(1-v)^2))}{(1-(1-v)^2)}\right)\left(\frac{(1-(1-v')^3) - (1-(1-v)^3)}{(1-(1-v)^3)}\right)\left(\frac{v}{v'-v}\right)\right)$$

**Example 2:**

Let's do example 1 numerically:

$$X_0 = (\mu_0, v_0) = (0.1, 0.7)$$
$$Y_0 = (\mu'_0, v'_0) = (0.6, 0.3)$$
$$X_0 \trianglelefteq Y_0 \Leftrightarrow X_0 \oplus \varepsilon = Y_0$$
$$\Leftrightarrow \varepsilon = Y_0 \ominus X_0$$
$$= (0.6, 0.3) \ominus (0.1, 0.7)$$
$$= \left(\frac{0.6 - 0.1}{1 - 0.1}, \frac{0.3}{0.7}\right) IFN$$
$$= \left(\frac{5}{9}, \frac{3}{7}\right) IFN$$

$$\varphi(X) = X^2 = (\mu, v)^2 = (\mu^2, 1 - (1-v)^2) = (f(\mu), g(v))$$

$$\varphi(Y) = Y^2 = (\mu', v')^2 = (\mu'^2, 1 - (1 - v')^2) = (f(\mu'), g(v'))$$

$$\gamma(X) = X^3 = (\mu, v)^3 = (\mu^3, 1 - (1 - v)^3) = (f_1(\mu), g_1(v))$$

$$\gamma(Y) = Y^3 = (\mu', v')^3 = (\mu'^3, 1 - (1 - v')^3) = (f_1(\mu'), g_1(v'))$$

$$\varphi(X_0) = X_0^2 = (0.1, 0.7)^2 = (0.01, 0.91) = (f(\mu), g(v))$$

$$\varphi(Y_0) = Y_0^2 = (0.6, 0.3)^2 = (0.36, 0.51) = (f(\mu'), g(v'))$$

$$\gamma(X_0) = X_0^3 = (0.1, 0.7)^3 = (0.001, 0.973) = (f_1(\mu), g_1(v))$$

$$\gamma(Y_0) = Y_0^3 = (0.6, 0.3)^3 = (0.216, 0.657) = (f_1(\mu'), g_1(v'))$$

**LHS:**

$$[\varphi(Y_0) \ominus \varphi(X_0)] \otimes \frac{d^\oplus \gamma(X)}{dX}$$

$$= \left( \frac{f_1(\mu_0') - f_1(\mu_0)}{1 - f_1(\mu_0)} \left( \frac{1 - \mu_0}{1 - f(\mu_0)} \frac{df(\mu)}{d\mu} \right), 1 + \left( \frac{g_1(v_0') - g_1(v_0)}{g_1(v_0)} \right) \left( \frac{v_0}{g(v_0)} \frac{dg(v)}{dv} \right) \right)$$

$$= \left( \frac{0.216 - 0.001}{1 - 0.001} \left( \frac{1 - 0.1}{1 - 0.01} \frac{df(\mu)}{d\mu} \right), 1 + \left( \frac{0.657 - 0.973}{0.973} \right) \left( \frac{0.7}{0.91} \frac{dg(v)}{dv} \right) \right)$$

$$= \left( \frac{215}{999} \left( \frac{90}{99} \frac{df(\mu)}{d\mu} \right), 1 - \left( \frac{316}{973} \right) \left( \frac{70}{91} \frac{dg(v)}{dv} \right) \right)$$

$$= \left( \frac{215}{999} \left( \frac{90}{99} \right) \frac{f(\mu') - f(\mu)}{\mu' - \mu}, 1 - \left( \frac{316}{973} \right) \left( \frac{70}{91} \frac{g(v') - g(v)}{v' - v} \right) \right)$$

$$= \left( \frac{215}{999} \left( \frac{90}{99} \right) \frac{0.36 - 0.01}{0.6 - 0.1}, 1 - \left( \frac{316}{973} \right) \left( \frac{70}{91} \frac{0.51 - 0.91}{0.3 - 0.7} \right) \right)$$

$$= (0.1369551369551, 0.7501778796743)$$

**RHS:**

$$(\gamma(Y) \ominus \gamma(X)) \otimes \frac{d^\oplus \varphi(X)}{dX}$$

$$= \left( \frac{f(\mu') - f(\mu)}{1 - f(\mu)} \left( \frac{1 - \mu}{1 - f_1(\mu)} \frac{df_1(\mu)}{d\mu} \right), 1 + \left( \frac{g(v') - g(v)}{g(v)} \right) \left( \frac{v}{g_1(v)} \frac{dg_1(v)}{dv} \right) \right)$$

$$= \left( \frac{0.36 - 0.01}{1 - 0.01} \left( \frac{1 - 0.1}{1 - 0.001} \frac{df_1(\mu)}{d\mu} \right), 1 + \left( \frac{0.51 - 0.91}{0.91} \right) \left( \frac{0.7}{0.973} \frac{dg_1(v)}{dv} \right) \right)$$

$$= \left(\frac{35}{99}\left(\frac{900}{999}\frac{df_1(\mu)}{d\mu}\right), 1 - \left(\frac{40}{91}\right)\left(\frac{700}{973}\frac{dg_1(v)}{dv}\right)\right)$$

$$= \left(\frac{35}{99}\left(\frac{900}{999}\left(\frac{f_1(\mu')-f_1(\mu)}{\mu'-\mu}\right)\right), 1 - \left(\frac{40}{91}\right)\left(\frac{700}{973}\left(\frac{g_1(v')-g_1(v)}{v'-v}\right)\right)\right)$$

$$= \left(\frac{35}{99}\left(\frac{900}{999}\left(\frac{0.216-0.001}{0.6-0.1}\right)\right), 1 - \left(\frac{40}{91}\right)\left(\frac{700}{973}\left(\frac{0.657-0.973}{0.3-0.7}\right)\right)\right)$$

$$= (0.1369551369551, 0.7501778796743)$$

### Corollary 3.2:

The mean value theorem is a special case of the cauchy mean value theorem, that is for $\gamma(X) = X$, then CMVT reduces to MVT.

### Proof:

Firstly, We should show that if $\gamma(X) = X = (\mu, v) = (f(\mu), g(v))$, then $\frac{d^\oplus \gamma(X)}{dX} = E = (1,0)$.

That is,

$$\frac{d^\oplus \gamma(X)}{dX} = \left(\frac{1-\mu}{1-f(\mu)}\frac{df(\mu)}{d\mu}, 1 - \frac{v}{g(v)}\frac{dg(v)}{dv}\right)$$

$$= \left(\frac{1-\mu}{1-\mu}\frac{df(\mu)}{d\mu}, 1 - \frac{v}{v}\frac{dg(v)}{dv}\right)$$

$$= \left(\frac{1-\mu}{1-\mu}\frac{d\mu}{d\mu}, 1 - \frac{v}{v}\frac{dv}{dv}\right)$$

$$= (1,0) = E$$

Let $\gamma(X) = X$,

Then we have:

$$[\varphi(Y) \ominus \varphi(X)] \otimes \frac{d^\oplus \gamma(X)}{dX} = [\gamma(Y) \ominus \gamma(X)] \otimes \frac{d^\oplus \varphi(X)}{dX}$$

$$\Rightarrow [\varphi(Y) \ominus \varphi(X)] \otimes E = (Y \ominus X) \otimes \frac{d^\oplus \varphi(X)}{dX}$$

$$\Leftrightarrow [\varphi(Y) \ominus \varphi(X)] = \frac{d^{\oplus}\varphi(X)}{dX} \otimes (Y \ominus X)$$

### Corollary 3.3:

The Rolle's theorem is a special case of Cauchy mean value theorem, that is for $\varphi(X) = \varphi(Y)$, $X \neq Y$, then CMVT reduces to Rolle's theorem.

**Proof:**

Let $\varphi(X) = \varphi(Y)$, then

$$(Y \ominus X) \otimes \frac{d^{\oplus}\varphi(X)}{dX} = [\varphi(Y) \ominus \varphi(X)] = (0.1) = 0$$

$$\Rightarrow \frac{d^{\oplus}\varphi(X)}{dX} = (0,1) = 0$$

QED.

**Example:**

$$X = (0.2, 0.4)$$
$$Y = (0.4, 0.2)$$
$$\varphi(X) = \alpha_0 = (f(\mu_0), g(v_0)) = Constant$$

$$\frac{d^{\oplus}\varphi(X)}{dX} = \left( \frac{1-\mu}{1-f(\mu_0)} \frac{df(\mu_0)}{d\mu}, 1 - \frac{v}{g(v)} \frac{dg(v_0)}{dv} \right)$$

$$= \left( \frac{1-\mu}{1-f(\mu_0)} 0, 1 - \frac{v}{g(v)} 0 \right)$$

$$= (0,1)$$

$$= \mathbf{0}$$

Using rolle's theorem:

$$\frac{d^{\oplus}\varphi(X)}{dX} = [\varphi(Y) \ominus \varphi(X)] = \alpha_0 \ominus \alpha_0 = (0.1)$$

### Corollary 3.5:

Let $\varphi(X) = (f(\mu), g(v))$ and $\gamma(X) = (f_1(\mu), g_1(v))$ be a monotone increasing IFF of $X \unlhd Y$. If $\varphi(X)$ and $\gamma(X)$ are differentiable and satisfy CMVT, then $\lambda\varphi(X)$ and $\lambda\gamma(X)$ are also satisfy CMVT.

**Proof:**

Firstly, we should show

$$\frac{d^{\oplus}(\lambda\varphi(X))}{dX} = (\lambda, 1-\lambda) \otimes \frac{d^{\oplus}\varphi(X)}{dX}$$

That is:

$$\lambda\varphi(X) = \lambda(f(\mu), g(v)) = (1 - (1-f(\mu))^{\lambda}, (g(v))^{\lambda})$$

$$\frac{d^{\oplus}\varphi(X)}{dX} = \left(\frac{1-\mu}{1-f(\mu)}\frac{df(\mu)}{d\mu}, 1 - \frac{v}{g(v)}\frac{dg(v)}{dv}\right)$$

$$\frac{d^{\oplus}\lambda\varphi(X)}{dX} = \left(\frac{1-\mu}{1-(1-(1-f(\mu))^{\lambda})}\frac{d(1-(1-f(\mu))^{\lambda})}{d\mu}, 1 - \frac{v}{(g(v))^{\lambda}}\frac{d((g(v))^{\lambda})}{dv}\right)$$

$$= \left(\frac{1-\mu}{(1-f(\mu))^{\lambda}}\lambda(1-f(\mu))^{\lambda-1}\frac{df(\mu)}{d\mu}, 1 - \frac{v}{(g(v))^{\lambda}}\lambda(g(v))^{\lambda-1}\frac{dg(v)}{dv}\right)$$

$$= \left(\lambda\frac{1-\mu}{1-f(\mu)}\frac{df(\mu)}{d\mu}, \left(1 - \lambda\frac{v}{g(v)}\frac{dg(v)}{dv}\right)\right)$$

and

$$(\lambda, 1-\lambda) \otimes \frac{d\gamma(X)}{dX} = (\lambda, 1-\lambda) \otimes \left( \frac{1-\mu}{1-f(\mu)} \frac{df(\mu)}{d\mu}, 1 - \frac{v}{g(v)} \frac{dg(v)}{dv} \right)$$

$$= (\lambda, 1-\lambda) \otimes \left( \lambda \frac{1-\mu}{1-f(\mu)} \frac{df(\mu)}{d\mu}, (1-\lambda) + \left(1 - \frac{v}{g(v)} \frac{dg(v)}{dv}\right) \right.$$
$$\left. - (1-\lambda)\left(1 - \frac{v}{g(v)} \frac{dg(v)}{dv}\right) \right)$$

$$= \left( \lambda \frac{1-\mu}{1-f(\mu)} \frac{df(\mu)}{d\mu}, \left(1 - \frac{v}{g(v)} \frac{dg(v)}{dv}\right) + (1-\lambda)\left(\frac{v}{g(v)} \frac{dg(v)}{dv}\right) \right)$$

$$= \left( \lambda \frac{1-\mu}{1-f(\mu)} \frac{df(\mu)}{d\mu}, \left(1 - \lambda \frac{v}{g(v)} \frac{dg(v)}{dv}\right) \right)$$

Then we have

$$\frac{d^{\oplus}\lambda\varphi(X)}{dX} = (\lambda, 1-\lambda) \otimes \frac{d^{\oplus}\varphi(X)}{dX}$$

Also

$$\frac{d^{\oplus}\lambda\gamma(X)}{dX} = (\lambda, 1-\lambda) \otimes \frac{d^{\oplus}\gamma(X)}{dX}$$

Let $\varphi(X)$ and $\gamma(X)$ satisfy CMVT, then we have:

$$[\varphi(Y) \ominus \varphi(X)] \otimes \frac{d^{\oplus}\gamma(X)}{dX} = (\gamma(Y) \ominus \gamma(X)) \otimes \frac{d^{\oplus}\varphi(X)}{dX}$$

$$\Leftrightarrow [\lambda(\varphi(Y) \ominus \varphi(X))] \otimes (\lambda, 1-\lambda) \otimes \frac{d^{\oplus}\gamma(X)}{dX} = (\lambda(\gamma(Y) \ominus \gamma(X))) \otimes (\lambda, 1-\lambda) \otimes \frac{d^{\oplus}\varphi(X)}{dX}$$

$$\Leftrightarrow [\lambda\varphi(Y) \ominus \lambda\varphi(X)] \otimes \frac{d^{\oplus}\lambda\gamma(X)}{dX} = (\lambda\gamma(Y) \ominus \lambda\gamma(X)) \otimes \frac{d^{\oplus}\lambda\varphi(X)}{dX}$$

So that, $\lambda\varphi(X)$ and $\lambda\gamma(X)$ also satisfy CMVT.

QED.

**Corollary 3.6:**

Let $\alpha$ and $\beta$ are IFNs, $\varphi(X) = (f(\mu), g(v))$ and $\gamma(X) = (f_1(\mu), g_1(v))$ be a monotone increasing IFFs of $X \trianglelefteq Y$. If $\varphi(X)$ and $\gamma(X)$ are differentiable and satisfy CMVT, then $\alpha \oplus \varphi(X)$ and $\beta \oplus \gamma(X)$ are also satisfy CMVT.

**Proof:**

Firstly we should show

$$\alpha \oplus \varphi(X) = (\mu_\alpha, v_\alpha) \oplus (f(\mu), g(v))$$
$$= (\mu_\alpha + f(\mu) - \mu_\alpha f(\mu), v_\alpha g(v))$$

$$\frac{d^\oplus(\alpha \oplus \varphi(X))}{dX} = \left( \frac{1-\mu}{1-(\mu_\alpha + f(\mu) - \mu_\alpha f(\mu))} \frac{d(\mu_\alpha + f(\mu) - \mu_\alpha f(\mu))}{d\mu}, 1 \right.$$
$$\left. - \frac{v}{(v_\alpha g(v))} \frac{d(v_\alpha g(v))}{dv} \right)$$

$$= \left( \frac{1-\mu}{(1-\mu_\alpha)(1-f(\mu))} (1-\mu_\alpha) \frac{d(f(\mu))}{d\mu}, 1 - \frac{v_\alpha v}{(v_\alpha g(v))} \frac{d(g(v))}{dv} \right)$$

$$= \left( \frac{1-\mu}{(1-f(\mu))} \frac{d(f(\mu))}{d\mu}, 1 - \frac{v}{g(v)} \frac{dg(v)}{dv} \right)$$

$$= \frac{d^\oplus \varphi(X)}{dX}$$

And also

$$\frac{d^\oplus(\beta \oplus \gamma(X))}{dX} = \frac{d^\oplus \varphi(X)}{dX}$$

Now,

$$[(\alpha \oplus \varphi(Y)) \ominus (\alpha \oplus \varphi(X))] \otimes \frac{d^\oplus(\beta \oplus \gamma(X))}{dX} = ((\beta \oplus \gamma(Y)) \ominus (\beta \oplus \gamma(X))) \otimes \frac{d^\oplus(\alpha \oplus \varphi(X))}{dX}$$

$$\Leftrightarrow [\varphi(Y) \ominus \varphi(X)] \otimes \frac{d^\oplus \gamma(X)}{dX} = (\gamma(Y) \ominus \gamma(X)) \otimes \frac{d^\oplus \varphi(X)}{dX}$$

QED

## 3.2 Multiplication CMVT

**Theorem 3.1** [MVT] Let $\varphi(X) = (f(\mu), g(v))$ be a monotone increasing IFF of $X \trianglelefteq Y$. If $\varphi(X)$ is differentiable, then:

$$\varphi(Y) \oslash \varphi(X) = \frac{d^{\otimes}\varphi(X)}{dX} \oplus (Y \oslash X)$$

**Proof:**

If $X \trianglelefteq Y$, $X$ and $Y$ are expressed by $(\mu, v)$ and $(\mu', v')$ respectively, then

$$\Delta X = Y \oslash X = (\mu', v') \oslash (\mu, v) = \left(\frac{\mu'}{\mu}, \frac{v'-v}{1-v}\right)$$

and also,

$$\frac{d^{\otimes}\varphi(X)}{dX} = \left(1 - \frac{\mu}{f(\mu)}\frac{df(\mu)}{d\mu}, \frac{1-v}{1-g(v)}\frac{dg(v)}{dv}\right)$$

Hence, we have:

$$\varphi(Y) \oslash \varphi(X) = \big(f(\mu'), g(v')\big) \oslash \big(f(\mu), g(v)\big)$$
$$= \left(\frac{f(\mu')}{f(\mu)}, \frac{g(v')-g(v)}{1-g(v)}\right)$$

And,

$$\frac{d^{\otimes}\varphi(X)}{dX} \oplus \Delta X = \left(1 - \frac{\mu}{f(\mu)}\frac{df(\mu)}{d\mu}, \frac{1-v}{1-g(v)}\frac{dg(v)}{dv}\right) \oplus \left(\frac{\mu'}{\mu}, \frac{v'-v}{1-v}\right)$$
$$= \left(1 - \frac{\mu}{f(\mu)}\frac{df(\mu)}{d\mu} + \frac{\mu'}{\mu} - \left(1 - \frac{\mu}{f(\mu)}\frac{df(\mu)}{d\mu}\right)\frac{\mu'}{\mu}, \frac{1-v}{1-g(v)}\frac{dg(v)}{dv}\frac{v'-v}{1-v}\right)$$
$$= \left(1 - \frac{\mu-\mu'}{f(\mu)}\frac{df(\mu)}{d\mu}, \frac{v'-v}{1-g(v)}\frac{dg(v)}{dv}\right) = \left(1 + \frac{f(\mu')-f(\mu)}{f(\mu)}, \frac{g(v')-g(v)}{1-g(v)}\right)$$
$$= \left(\frac{f(\mu')}{f(\mu)}, \frac{g(v')-g(v)}{1-g(v)}\right) = \varphi(Y) \oslash \varphi(X)$$

Therefore:

$$\varphi(Y) \oslash \varphi(X) = \frac{d^\otimes \varphi(X)}{dX} \oplus (Y \oslash X)$$

QED.

**Example:** Show that Multiplication MVT can be written as

$$\frac{d^\otimes \varphi(X)}{dX} = \frac{\varphi(Y)}{\varphi(X)} \ominus \frac{Y}{X}$$

Solution:

$$\frac{\varphi(Y)}{\varphi(X)} \ominus \frac{Y}{X} = (\varphi(Y) \oslash \varphi(X)) \ominus (Y \oslash X) = \left(\frac{f(\mu')}{f(\mu)}, \frac{g(v') - g(v)}{1 - g(v)}\right) \ominus \left(\frac{\mu'}{\mu}, \frac{v' - v}{1 - v}\right)$$

$$= \left(\frac{\frac{f(\mu')}{f(\mu)} - \frac{\mu'}{\mu}}{1 - \frac{\mu'}{\mu}}, \frac{g(v') - g(v)}{1 - g(v)} \frac{1 - v}{v' - v}\right)$$

$$= \left(-\frac{\frac{\mu'}{\mu} - \frac{f(\mu')}{f(\mu)}}{1 - \frac{\mu'}{\mu}}, \frac{1 - v}{1 - g(v)} \frac{g(v') - g(v)}{v' - v}\right)$$

$$= \left(\frac{\mu' f(\mu) - \mu f(\mu')}{f(\mu)\mu' - \mu}, \frac{1 - v}{1 - g(v)} \frac{g(v') - g(v)}{v' - v}\right)$$

$$= \left(\frac{\mu' f(\mu) - \mu f(\mu')}{f(\mu)(\mu' - \mu)}, \frac{1 - v}{1 - g(v)} \frac{g(v') - g(v)}{v' - v}\right)$$

$$= \left(\frac{\mu' f(\mu) - \mu f(\mu')}{f(\mu)(\mu' - \mu)}, \frac{1 - v}{1 - g(v)} \frac{g(v') - g(v)}{v' - v}\right)$$

$$
\begin{aligned}
&= \left( \frac{\mu' f(\mu) - \mu f(\mu) - \mu f(\mu') + \mu f(\mu)}{f(\mu)(\mu' - \mu)}, \frac{1-v}{1-g(v)} \frac{dg(v)}{dv} \right) \\
&= \left( \frac{f(\mu)(\mu' - \mu) - \mu f(\mu') + \mu f(\mu)}{f(\mu)(\mu' - \mu)}, \frac{1-v}{1-g(v)} \frac{dg(v)}{dv} \right) \\
&= \left( \frac{f(\mu)(\mu' - \mu)}{f(\mu)(\mu' - \mu)} - \frac{\mu f(\mu') - \mu f(\mu)}{f(\mu)(\mu' - \mu)}, \frac{1-v}{1-g(v)} \frac{dg(v)}{dv} \right) \\
&= \left( 1 - \frac{\mu f(\mu') - \mu f(\mu)}{f(\mu)(\mu' - \mu)}, \frac{1-v}{1-g(v)} \frac{dg(v)}{dv} \right) \\
&= \left( 1 - \frac{\mu}{f(\mu)} \frac{f(\mu') - f(\mu)}{\mu' - \mu}, \frac{1-v}{1-g(v)} \frac{dg(v)}{dv} \right) \\
&= \left( 1 - \frac{\mu}{f(\mu)} \frac{df(\mu)}{d\mu}, \frac{1-v}{1-g(v)} \frac{dg(v)}{dv} \right) = \frac{d^{\otimes}\varphi(X)}{dX}
\end{aligned}
$$

Based on theorem 3.3 we generalize MVT into CMVT, as follows:

**Theorem 3.4**

Let $\varphi(X) = (f(\mu), g(v))$ and $\gamma(X) = (f_1(\mu), g_1(v))$ be a monotone increasing IFF of $X \trianglelefteq Y$. If $\varphi(X)$ and $\gamma(X)$ are differentiable, then

$$[\varphi(Y) \oslash \varphi(X)] \oplus \frac{d^{\otimes}\gamma(X)}{dX} = \frac{d^{\otimes}\varphi(X)}{dX} \oplus [\gamma(Y) \oslash \gamma(X)]$$

**Proof:**

Let:
$$X = (\mu, v)$$
$$Y = (\mu', v')$$
$$\varphi(X) = (f(\mu), g(v))$$
$$\varphi(Y) = (f(\mu'), g(v'))$$
$$\gamma(X) = (f_1(\mu), g_1(v))$$
$$\gamma(Y) = (f_1(\mu'), g_1(v'))$$

$$\frac{d^{\otimes}\varphi(X)}{dX} = \left(1 - \frac{\mu}{f(\mu)}\frac{df(\mu)}{d\mu}, \frac{1-v}{1-g(v)}\frac{dg(v)}{dv}\right)$$

$$\frac{d^{\otimes}\gamma(X)}{dX} = \left(1 - \frac{\mu}{f_1(\mu)}\frac{df_1(\mu)}{d\mu}, \frac{1-v}{1-g_1(v)}\frac{dg_1(v)}{dv}\right)$$

Then we have:

$$\varphi(Y) \oslash \varphi(X) = \big(f(\mu'), g(v')\big) \oslash \big(f(\mu), g(v)\big)$$

$$= \left(\frac{f(\mu')}{f(\mu)}, \frac{g(v') - g(v)}{1 - g(v)}\right)$$

And also:

$$\gamma(Y) \oslash \gamma(X)$$
$$= \big(f_1(\mu'), g_1(v')\big) \oslash \big(f_1(\mu), g_1(v)\big)$$
$$= \left(\frac{f_1(\mu')}{f_1(\mu)}, \frac{g_1(v') - g_1(v)}{1 - g_1(v)}\right)$$

So that, the left hand side is:

$$[\varphi(Y) \oslash \varphi(X)] \oplus \frac{d^{\otimes}\gamma(X)}{dX}$$

$$= \left(\frac{f(\mu')}{f(\mu)}, \frac{g(v') - g(v)}{1 - g(v)}\right) \oplus \left(1 - \frac{\mu}{f_1(\mu)}\frac{df_1(\mu)}{d\mu}, \frac{1-v}{1-g_1(v)}\frac{dg_1(v)}{dv}\right)$$

$$= \left(\frac{f(\mu')}{f(\mu)} + 1 - \frac{\mu}{f_1(\mu)}\frac{df_1(\mu)}{d\mu}\right.$$

$$\left. - \frac{f(\mu')}{f(\mu)}\left(1 - \frac{\mu}{f_1(\mu)}\frac{df_1(\mu)}{d\mu}\right), \left(\frac{g(v') - g(v)}{1 - g(v)}\right)\left(\frac{1-v}{1-g_1(v)}\frac{dg_1(v)}{dv}\right)\right)$$

$$= \left(\frac{f(\mu')}{f(\mu)} + 1 - \frac{\mu}{f_1(\mu)}\frac{df_1(\mu)}{d\mu} - \frac{f(\mu')}{f(\mu)}\right.$$

$$\left. + \frac{f(\mu')}{f(\mu)}\left(\frac{\mu}{f_1(\mu)}\frac{df_1(\mu)}{d\mu}\right), \left(\frac{g(v') - g(v)}{1 - g(v)}\right)\left(\frac{1-v}{1-g_1(v)}\frac{dg_1(v)}{dv}\right)\right)$$

$$= \left(1 - \frac{\mu}{f_1(\mu)}\frac{df_1(\mu)}{d\mu} + \frac{f(\mu')}{f(\mu)}\left(\frac{\mu}{f_1(\mu)}\frac{df_1(\mu)}{d\mu}\right), \left(\frac{g(v') - g(v)}{1 - g(v)}\right)\left(\frac{1-v}{1-g_1(v)}\frac{dg_1(v)}{dv}\right)\right)$$

$$= \left(1 + \left(\frac{f(\mu') - f(\mu)}{f(\mu)}\right)\left(\frac{\mu}{f_1(\mu)}\frac{df_1(\mu)}{d\mu}\right), \left(\frac{g(v') - g(v)}{1 - g(v)}\right)\left(\frac{1-v}{1-g_1(v)}\frac{dg_1(v)}{dv}\right)\right)$$

$$= \left(1 + \left(\frac{f(\mu') - f(\mu)}{f(\mu)}\right)\frac{\mu}{f_1(\mu)}\left(\frac{f_1(\mu') - f_1(\mu)}{\mu' - \mu}\right),\right.$$

$$\left.\left(\frac{g(v') - g(v)}{1 - g(v)}\right)\left(\frac{1-v}{1-g_1(v)}\frac{g_1(v') - g_1(v)}{v' - v}\right)\right)$$

$$= \left(1 + \frac{f(\mu') - f(\mu)}{f(\mu)}\frac{f_1(\mu') - f_1(\mu)}{f_1(\mu)}\left(\frac{\mu}{\mu' - \mu}\right), \frac{g(v') - g(v)}{1 - g(v)}\frac{g_1(v') - g_1(v)}{1 - g_1(v)}\left(\frac{1-v}{v' - v}\right)\right)$$

The right hand side is:

$$[\gamma(Y) \oslash \gamma(X)] \oplus \frac{d^\otimes \varphi(X)}{dX}$$

$$= \left(\frac{f_1(\mu')}{f_1(\mu)}, \frac{g_1(v') - g_1(v)}{1 - g_1(v)}\right) \oplus \left(1 - \frac{\mu}{f(\mu)}\frac{df(\mu)}{d\mu}, \frac{1-v}{1-g(v)}\frac{dg(v)}{dv}\right)$$

$$= \left(\frac{f_1(\mu')}{f_1(\mu)} + 1 - \frac{\mu}{f(\mu)}\frac{df(\mu)}{d\mu} - \frac{f_1(\mu')}{f_1(\mu)}\left(1 - \frac{\mu}{f(\mu)}\frac{df(\mu)}{d\mu}\right), \frac{g_1(v') - g_1(v)}{1 - g_1(v)}\left(\frac{1-v}{1-g(v)}\frac{dg(v)}{dv}\right)\right)$$

$$= \left(1 + \frac{f_1(\mu') - f_1(\mu)}{f_1(\mu)}\frac{\mu}{f(\mu)}\frac{df(\mu)}{d\mu}, \frac{g_1(v') - g_1(v)}{1 - g_1(v)}\left(\frac{1-v}{1-g(v)}\frac{dg(v)}{dv}\right)\right)$$

$$= \left(1 + \frac{f_1(\mu') - f_1(\mu)}{f_1(\mu)}\frac{\mu}{f(\mu)}\frac{f(\mu') - f(\mu)}{\mu' - \mu}, \frac{g_1(v') - g_1(v)}{1 - g_1(v)}\left(\frac{1-v}{1-g(v)}\frac{g(v') - g(v)}{v' - v}\right)\right)$$

$$= \left(1 + \frac{f_1(\mu') - f_1(\mu)}{f_1(\mu)}\frac{f(\mu') - f(\mu)}{f(\mu)}\frac{\mu}{\mu' - \mu}, \frac{g_1(v') - g_1(v)}{1 - g_1(v)}\frac{g(v') - g(v)}{1 - g(v)}\left(\frac{1-v}{v' - v}\right)\right)$$

Now, the left hand side and the right hand side are equal.

**QED**

## 4. Potential Application In Decision Analysis

In big data era, decision making needs an ability to analyze the changes of trends that run rapidly and contain complicated information. In intuitionistic fuzzy calculus IFC with domain IFN, then the roles become will have some benefits. Where the situation faced by decision maker is changes follow Intuitionistic Fuzzy Function IFF formulation, then we need to know how the influence or relation of changes IFN variable to IFF. We can use Intuitionistic Fuzzy Derivative IFD to answer the problem.

The decision making's flow chart can be seen bellow:

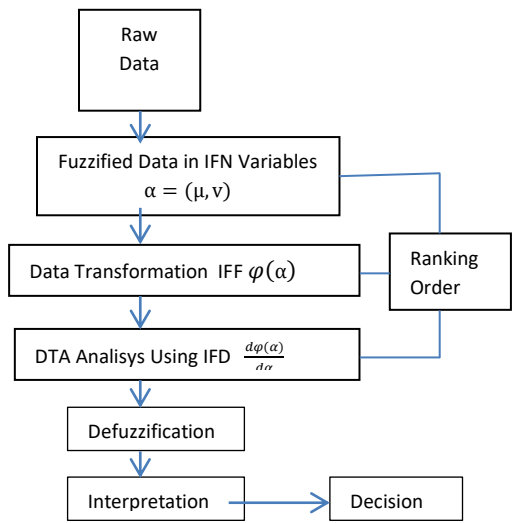

**Fig. 4.1 Flowchart Decision Trend Analysis (DTA) Using IFC**

The CMVT has a nice property about monotonically IFFs, so that DTA may find advantages in doing analysis on decision trend. Even tough, there may still along journey to find practical application on DTA.

## 5. Conclusion

In this paper, we have first investigated the Xu monotone IFF theorem. In the theorem mean value was written as approximation equation rather than in an equality. Based on it, we have derived a mean value theorem in intuitionistic fuzzy calculus environment in the form of Addition MVT and Multiplication MVT. We also derived and proved the general version of the mean value theorem of IFF in Cauchy version and we have proved two theorems, Addition CMVT and Multiplication CMVT. We give some examples for the theorems to show the calculation process. We have found that from the CMVT of IFF emerged some interesting corollaries. The most important corollary is that we can derive MVT and Rolle's theorem from CMVT. Other nice properties also have been shown that some operations are invariant on CMVT, that is scalar multiplication and IFN addition on IFFs also satisfy CMVT. Finally, we consider the potential application on Decision trend analysis using intuitionistic fuzzy variable.

## 6. Acknowledgement

This article will be part of Phd dissertation at Department of Mathematics and Data Sciences of Andalas University. I must thank also the anonymous reviewer for his/her valuable comments and suggestions which inspired the results of Section 3.